\documentclass{article}
\usepackage{amsmath, amscd, amsfonts, amssymb, amsthm}

\newtheorem{theorem}{Theorem}[section]
\newtheorem{corollary}[theorem]{Corollary}
\newtheorem{lemma}[theorem]{Lemma}
\newtheorem{prop}[theorem]{Proposition}
\begin{document}

\title{Perturbations of the Symmetric Exclusion Process\footnote{Research supported by NSF grant DMS-00-70465}}
\author{Paul Jung}

\date{\textit{Department of Mathematics, Cornell University, Ithaca, NY 14853, USA}} \maketitle

\begin{abstract}
This paper gives results concerning the asymptotics of the
invariant measures, $\mathcal{I},$ for exclusion processes where
$p(x,y)=p(y,x)$ except for finitely many $x,y\in\mathcal{S}$ and
$p(x,y)$ corresponds to a transient Markov chain on $\mathcal{S}$.
As a consequence, a complete characterization of $\mathcal{I}$ is
given for the case where $p(x,y)=p(y,x)$ for all but a single
ordered pair $(u,v)$. Also, this paper addresses the question:
When do local changes to a symmetric kernel $p(x,y)=p(y,x)$ affect
the evolution of the exclusion process globally?
\end{abstract}


\section{Introduction}
The exclusion process is a well-known interacting particle system
that has been used in biology as a model for the particle motion
of ribosomes (Macdonald, Gibbs, and Pipkin(1968)), in physics as a
model for a lattice gas at infinite temperature (Spitzer(1970)),
and in ecology as a model in which two opposing species swap
territory (Clifford and Sudbury(1973)). The state space for the
exclusion process is $X=\{0,1\}^{\mathcal{S}}$ for $\mathcal{S}$ a
countable set, and its generator is given by the closure of the
operator $\Omega$ on $\mathcal{D}(X)$, the set of all functions on
$X$ depending on finitely many coordinates. Let
\begin{equation*}
\sup_y \sum_x p(x,y)<\infty \text{ and }\sup_x \sum_y
p(x,y)<\infty\text{ for }p(x,y)\ge 0.
\end{equation*}
If $f\in \mathcal{D}(X)$ and
\begin{equation*}
\eta_{xy}(u)=\left\{
\begin{array}{ll}
\eta(y)&\text{if }u=x\\
\eta(x)&\text{if }u=y\\
\eta(u)&\text{if }u\neq x,y\\
\end{array}
\right.
\end{equation*}
then
\begin{equation}\label{qgenerator}
{\Omega}f(\eta)=\sum_{x,y}{p(x,y)\eta(x)(1-\eta(y))[f(\eta_{xy})-f(\eta)]}.
\end{equation}
We will denote the semigroup of this process by ${S}(t)$.

An intuitive description of the process is given by thinking of
the $1$'s as particles and the $0$'s as empty sites. A particle at
site $x\in \mathcal{S}$ waits an exponential time with parameter
$p(x)=\sum_y p(x,y)$ at which time it chooses a $y\in \mathcal{S}$
with probability $p(x,y)/p(x)$. If $y$ is empty then the particle
at $x$ goes to $y$, while if $y$ is occupied the particle at $x$
does not move.

The construction of the exclusion process is fully described in
IPS (Liggett(1985)). It is assumed there that the transition
kernel satisfies $\sum_y p(x,y)=1$, however, this is just a
normalization of the process we have just described.  To see this,
simply add self-jump rates to the process we have described above:
\begin{equation*}
p(x,x)= \sup_z \sum_y p(z,y)-\sum_y p(x,y).
\end{equation*}
Dividing all transition rates by $\sup_z \sum_y p(z,y)$ gives us
the process constructed in IPS.

Let $\nu_\alpha$ be the product measure on
$X=\{0,1\}^{\mathcal{S}}$ with marginals
$\nu_\alpha\{\eta:\eta(x)=1, x\in\mathcal{S}\}=\alpha(x)$. When
the transition kernel is irreducible and symmetric,
$p(x,y)=p(y,x)$ for all $x,y\in\mathcal{S}$, the set of extremal
invariant measures for the process is given by
\begin{equation}\label{q15}
\mathcal{I}_e=\{\lim_{t\rightarrow\infty}\nu_\alpha S(t): 0\le
\alpha(x)\le 1\text{ and }\sum_y p(x,y)\alpha(y)=\alpha(x)\text{
for all }x\in\mathcal{S}\}.
\end{equation}
The above characterization of $\mathcal{I}_e$ for symmetric
processes is carried out by studying the finite-particle exclusion
process which is the dual process of the infinite-particle
exclusion process. In fact, the limit of $\nu_\alpha S(t)$ as $t$
goes to infinity is known to exist because of this duality. One
should note that by the Krein-Milman theorem, characterizing
$\mathcal{I}_e$ is equivalent to characterizing $\mathcal{I}$. For
details on the symmetric exclusion process we refer the reader to
Chapter VIII of IPS.

If the transition kernel is not symmetric then the dual is not
available, and the problem of classifying $\mathcal{I}$ becomes
exceedingly more difficult.  In fact there are only a few cases
for which $\mathcal{I}$ is totally known.  We refer the reader to
Jung(2003) for a synopsis of those cases.

In this paper we will consider exclusion processes which are
perturbations of symmetric exclusion.  A perturbation of an
exclusion process with irreducible transition kernel $p(x,y)$ is
an exclusion process with transition kernel $\bar{p}(x,y)$
satisfying the following. Let $\bar{p}(x,y)=p(x,y)$ for all
$(x,y)$ except for $n$ ordered pairs
$\{(x_1,y_1),\ldots,(x_n,y_n)\}$. At $(x_i, y_i)$ we have the
perturbation $\bar{p}(x_i,y_i)=p(x_i,y_i)+\epsilon_i$ for
$\epsilon_i> -p(x_i,y_i)$.  Note that this implies
$\bar{p}(x,y)>0$ if and only if $p(x,y)>0$. Also, note that the
$x_i$'s and $y_i$'s are not necessarily distinct. If the original
kernel $p(x,y)$ is symmetric then we will say $\bar{p}(x,y)$ is
\textit{quasi-symmetric}.
%In order to avoid complications we
%will also assume hereafter that $\bar{p}(x,y)$ is irreducible.
Throughout the rest of the paper ${S}(t)$ and $\mathcal{I}$ will
denote the semigroup  and invariant measures corresponding to
$p(x,y)$ while $\bar{S}(t)$ and $\bar{\mathcal{I}}$ will be the
semigroup and invariant measures corresponding to the perturbed
kernel $\bar{p}(x,y)$.

As noted earlier, an analog of the dual finite-particle exclusion
process of the symmetric exclusion process in Chapter VIII of IPS
does not exist for quasi-symmetric processes which are not
symmetric. However, an approximation to the dual is available
which makes the the study of quasi-symmetric processes much more
tenable than processes with no symmetry whatsoever.  Also, the
fact that quasi-symmetric kernels are mostly symmetric allows us
to use a coupling technique to prove a convergence result as well
as a complete characterization of $\bar{\mathcal{I}}$ for certain
quasi-symmetric exclusion processes.

Let $\mathcal{S}_k$ to be the set of all subsets of $\mathcal{S}$
containing $k$ elements.

\begin{theorem}\label{qtheorem2}
Suppose $\bar{p}(x,y)$ is a quasi-symmetric, irreducible
transition kernel corresponding to a transient Markov chain, and
suppose that $\mathcal{I}$ is the set of invariant measures
corresponding to the symmetric kernel $p(x,y)$. Then

(a) for each $\bar{\mu}\in\bar{\mathcal{I}}$ there exists a
measure $\mu\in\mathcal{I}$ such that
\begin{equation}\label{q11}
\lim_{n\rightarrow\infty} |\bar{\mu}\{\eta(x)=1\text{ for all
}x\in A^{n}\}-\mu\{\eta(x)=1\text{ for all }x\in A^{n}\}|=0
\end{equation}
for all $k$ and all sequences $\{A^n\}, A^n\in\mathcal{S}_k$ such
that each $x\in\mathcal{S}$ is in finitely many $A^n$, and

(b)  for each $\mu\in{\mathcal{I}}$ there exists a measure
$\bar{\mu}\in\bar{\mathcal{I}}$ satisfying (\ref{q11}).
\end{theorem}

Since we have a characterization of $\mathcal{I}$ given by
(\ref{q15}), the measure $\mu\in\mathcal{I}$ in part (a) must be
unique.  If one could somehow show that
$\bar{\mu}\in\bar{\mathcal{I}}$ in part (b) is unique as well,
then we would have a one-to-one correspondence between
$\mathcal{I}$ and $\bar{\mathcal{I}}$ thereby giving us a
characterization of $\bar{\mathcal{I}}$.  In Theorem
\ref{newtheorem} we prove exactly this for quasi-symmetric kernels
that are symmetric except for one ordered pair $(u,v)$.

From the point of view of practicality, Theorem \ref{qtheorem2}
gives us as good of a characterization of $\bar{\mathcal{I}}$ as
one could could hope for.  The reason for this is that even if one
were to show that $\bar{\mu}$ in part (b) is unique for all
quasi-symmetric kernels, one would not expect to be able to
calculate
\begin{equation}\label{teq1}
\bar{\mu}\{\eta(x)=1\text{ for all }x\in A\}
\end{equation}
explicitly for each finite $A\subset \mathcal{S}$.  The best one
could hope for is to know the asymptotics of (\ref{teq1}) for some
sequence $\{A^n\}$ in $\mathcal{S}_k$. But Theorem \ref{qtheorem2}
already gives us this.

\begin{theorem}\label{qtheorem}
Under the hypotheses of Theorem \ref{qtheorem2} we have that
$$\lim_{t\rightarrow\infty} \mu
\bar{S}(t)=\bar{\mu}\in\bar{\mathcal{I}}$$ exists for each
$\mu\in\mathcal{I}$ and $\bar{\mu}$ satisfies (\ref{q11}).
\end{theorem}

\begin{theorem}\label{newtheorem}
In addition to the hypotheses of Theorem \ref{qtheorem2} assume
that $\bar{p}(x,y)=p(x,y)$ except for exactly one ordered pair
$(u,v)$.
% where $\bar{p}(u,v)=p(u,v)+\epsilon, \epsilon>-p(u,v)$.
Then
$$\bar{\mathcal{I}}=\{\lim_{t\rightarrow\infty} \mu
\bar{S}(t):\mu\in\mathcal{I}\}.$$
%and
%\begin{eqnarray*}
%\bar{\mu}\{\eta(x)=1\text{ for all }x\in
%A\}&=&\mu\{\eta(x)=1\text{ for all }x\in
%A\}\\
%&+&\epsilon\mu\{\eta(u)=1,\eta(v)=0\}\sum_{x\in A}
%[G(u,x)-G(v,x)].
%\end{eqnarray*}
\end{theorem}

%By the characterization of $\mathcal{I}$ given in (\ref{q15})
%together with part (ii) of the theorem,
%$$\lim_{t\rightarrow\infty} \mu_1
%\bar{S}(t)=\lim_{t\rightarrow\infty} \mu_2 \bar{S}(t)$$ if and
%only if $\mu_1=\mu_2$ for $\mu_i\in\mathcal{I}$. Thus for each
%$\mu\in\mathcal{I}$ we get a distinct limiting measure in the set
%$\bar{\mathcal{I}}=\{\lim_{t\rightarrow\infty} \mu
%\bar{S}(t):\mu\in\mathcal{I}\}$.

% Also, note that the Markov chain
%corresponding to $\bar{p}(x,y)$ is transient if and only if the
%Markov chain with respect to ${p}(x,y)$ is transient.

Theorem \ref{qtheorem} has an interesting consequence motivated by
the following question: Does a local perturbation of the dynamics
of a process have global consequences on the evolution?

The answer is affirmative for quasi-symmetric exclusion processes
with nearest-neighbor kernels on $\mathcal{S}=\mathbb{Z}$. To see
this, consider the simple case where
\begin{eqnarray*}
\bar{p}(x,y)=1/2 \text{ for all }(x,y)\neq (0,1) \text{ and
}\bar{p}(0,1)=1/2+\epsilon, \epsilon>0.
\end{eqnarray*}
Then we can use Theorem 1.1 of Jung(2003) to find that the only
extremal invariant measures are the product measures $\{\nu^c:
0\le c \le\infty\}$ with marginals
\begin{equation*}
\nu^c\{\eta:\eta(x)=1\}=\left\{
\begin{array}{ll}
\frac{c}{1+c}&\text{for }x\le 0\\
\\
\frac{c+2c\epsilon}{1+c+2c\epsilon}&\text{for }x>0.
\end{array}\right.
\end{equation*}
Let $\nu_\rho$ be the product measure with marginals
$\nu_\rho\{\eta:\eta(x)=1\}=\rho\in[0,1]$. If we choose a sequence
of times $\{T_n\}$ going to infinity so that
\begin{equation*}\label{q.5}
\lim_{n\rightarrow\infty}\frac{1}{T_n}\int_0^{T_n} \nu_{\rho}
\Bar{S}(t) dt=\mu_\rho
\end{equation*}
exists, then Proposition I.1.8 in IPS tells us that $\mu_\rho$ is
invariant. Therefore it must be a mixture  of the measures
$\{\nu^c: 0\le c\le \infty\}$. Consequently
$$\lim_{x\rightarrow\infty}\mu_\rho\{\eta:\eta(x)=1\}>\lim_{x\rightarrow
-\infty}\mu_\rho\{\eta:\eta(x)=1\},$$ however, this clearly shows
that the perturbation at the origin affects the evolution of the
process globally.

On the other hand, Theorem \ref{qtheorem} tells us that
$\lim_{t\rightarrow\infty}\mu \Bar{S}(t)$ is not very different
from $\mu\in\mathcal{I}$ when a quasi-symmetric kernel
$\bar{p}(x,y)$ corresponds to a transient Markov chain. Thus we
have a negative answer to the above question. Our final theorem
gives us an indication as to what might be a good criterion for
determining when local perturbations of recurrent, symmetric
random walk kernels on $\mathbb{Z}$ and $\mathbb{Z}^2$ can have
global effects.

Given a kernel $p(x,y)$, the canonical graph associated to the set
$\mathcal{S}$ is the graph created by drawing an edge between $x$
and $y$ whenever $p(x,y)>0$. The graph is said to be
\textit{transitive} if the automorphism group acts transitively on
the vertex set $\mathcal{S}$. Hereafter we assume that any
transitive graph associated to $\mathcal{S}$ and $p(x,y)$ is
endowed with a metric which is also transitive with respect to the
automorphism group.  Fixing one vertex to be the origin $o$, let
$d(A,o)=\inf_{x\in A}|x|$ where $|x|$ denotes the distance from
$x$ to $o$. We also assume that for any $N>0$, there exists $x$
such that $|x|>N$.

\begin{theorem}\label{newtheorem2}
Let $p(x,y)$ be a recurrent, symmetric kernel which gives rise to
a transitive graph on the set $\mathcal{S}$, and let
$\bar{p}(x,y)$ be a corresponding quasi-symmetric kernel. Suppose
$G_n(x,y)$ is the Green's function corresponding to $p(x,y)$ (the
expected number of visits from $x$ to $y$ in $n$ steps). If for
each $y\in\mathcal{S}$
$$\lim_{|x|\rightarrow\infty}[a(x+y)-a(x)]=0$$ where
$$a(x)=\lim_{n\rightarrow\infty}[G_n(o,o)-G_n(o,x)],$$
then for each $\bar{\mu}\in\bar{\mathcal{I}}$ there exists
$\mu\in\mathcal{I}$ such that
\begin{equation}\label{q11.1}
\lim_{n\rightarrow\infty} |\bar{\mu}\{\eta(x)=1\text{ for all
}x\in A^{n}\}-\mu\{\eta(x)=1\text{ for all }x\in A^{n}\}|=0
\end{equation}
for all $k$ and all sequences $\{A^n\}, A^n\in\mathcal{S}_k$ such
that $\lim_{n\rightarrow\infty} d(A^n,o)=\infty$.
\end{theorem}

The organization of the rest of the paper is as follows.  We will
prove Theorem 1.1 in Section 2 by approximating the dual process
of symmetric exclusion.  In Section 3 we introduce a coupling
technique which is key in the proofs of Theorems 1.2, 1.3, and
1.4. In the last section we will prove Theorems 1.2-1.4.

\section{Approximating the Dual Process}\label{4}
In this section we will use an approximation to the dual of
symmetric exclusion in order to prove Theorem \ref{qtheorem2}. We
will assume in this section that $\lim_{t\rightarrow\infty} \mu
\bar{S}(t)$ exists for $\mu\in\mathcal{I}$ as this will be proved
in last section.

We now describe the dual finite-particle system $A_t$ used in the
analysis of symmetric systems.  The process $A_t$ is just the
normal exclusion process with the added condition that its initial
state $A_0$ has finitely many sites where $\eta(x)=1$. We write
$|A_t|=n$ to denote the number of sites that are $1$'s. In
particular $A_t$ is a countable-state Markov chain that acts like
$n$ independent particles having transition rates $p(x,y)$, except
that when a particle tries to move to an occupied site its motion
is suppressed.

In the sequel, we will need to think of the exclusion process in a
different way so that we can couple $\eta_t$ and $A_t$. Using a
symmetric transition kernel, assign to the subset
$\{x,y\}\in\mathcal{S}_2$ an exponential clock with rate $p(x,y)$.
Since $p(x,y)=p(y,x)$, this assignment is well-defined. Each time
the exponential clock for $\{x,y\}$ goes off, the values for
$\eta(x)$ and $\eta(y)$ will switch. This motion describes the
symmetric exclusion process.

We can now couple $A_t$ with $\eta_t$ using this new description.
The process $A_t$ is equal to $A_0$ until the first time that an
exponential clock for $\{x,y\}$ with $x\in A_0$ and $y\notin A_0$
goes off. At that time $A_t$ becomes $(A_0\backslash x)\cup y$.
Let $A^T_t$ be the dual process running backwards in time starting
from time $T$ so that $A^T_t=A_{T-t}$. Since the exponential times
for $\{x,y\}$ are uniformly distributed on $[0,T]$, we can use the
same clocks for both $A_t$ and $A_t^T$. We then have that
\begin{equation}\label{q12}
\{\eta_T(x)=1\text{ for all }x\in A_0^T\}=\{\eta_0(x)=1\text{ for
all }x\in A_T^T\}.
\end{equation}
The informed reader may recognize the similarity between
(\ref{q12}) and Theorem VIII.1.1 in IPS (duality of the exclusion
process).

Notice that when $\eta(x)=\eta(y)=1$, switching values is the same
as not switching values.  For the symmetric exclusion process, we
can reinterpret this statement in the following way. When a
particle tries to move to an occupied site, instead of its motion
being suppressed, the two particles switch places. This idea gives
us:

\begin{prop}\label{qlem5}
Suppose $\{A^n\}$ is a sequence in $\mathcal{S}_k$. If each
$x\in\mathcal{S}$ belongs to finitely many $A^n$ and the symmetric
kernel $p(x,y)$ corresponds to a transient Markov chain on
$\mathcal{S}$, then for each fixed $z\in\mathcal{S}$
$$\lim_{n\rightarrow\infty} P^{A^n}(z\in A^n_t\text{ for some }t\ge
0)=0.$$
\end{prop}
\begin{proof}
Let $Z_1(t),\ldots, Z_k(t)$ be $k$ particles each following the
motions of a Markov chain on $\mathcal{S}$ with transition rates
$p(x,y)$.  If $Z_i(t)=x$ and $Z_j(t)=y$ then since
$p(x,y)=p(y,x)$, we can couple the two processes so that $Z_i(t)$
goes to $y$ at the same time that $Z_j(t)$ goes to $x$.  If
$A^n=\{Z_1^n(0),\ldots, Z_k^n(0)\}$, then using this coupling
$A^n_t=\{Z_1^n(t),\ldots, Z_k^n(t)\}$. Therefore
$$\lim_{n\rightarrow\infty}P^{A^n}(z\in A^n_t\text{ for some }t\ge
0)\le \lim_{n\rightarrow\infty}\sum_{i=1}^k
P^{Z_i^n(0)}(Z_i^n(t)=z\text{ for some }t\ge 0)=0.$$
\end{proof}
Let
\begin{equation*}
B=\{x\in\mathcal{S}: \bar{p}(x,y)\neq {p}(x,y) \text{ or
}\bar{p}(y,x)\neq {p}(y,x)\text{ for some }y\in\mathcal{S}\}.
\end{equation*}

We will now describe a process $\bar{A}_t$ which approximates the
process $A_t$. In order to make the process $\bar{A}_t$ Markovian,
we have to assume that the filtration of $\bar{A}_t(\bar{\omega})$
takes into account the path space of the quasi-symmetric process
$\eta_t[\bar{\omega}]$ with sample path $\bar{\omega}$ (we use
$[\cdot ]$ here since $(\cdot )$ has been reserved for
$x\in\mathcal{S}$). In other words, $\bar{A}_t(\bar{\omega})$ and
$\eta_t[\bar{\omega}]$ share the same probability space.  In
particular, $\bar{P}_\nu$ is the measure on the path space of the
quasi-symmetric process $\eta_t[\bar{\omega}]$ having $\nu$ as its
initial distribution (likewise ${P}_\nu$ is the measure on the
path space of the symmetric process $\eta_t[\omega]$ with sample
path $\omega$). In order to avoid unnecessary technicalities we
will assume that the initial state $\bar{A}_0$ satisfies
$\bar{A}_0\cap B=\emptyset$.

We now describe the transitions of the process $\bar{A}_t$. If
$x\in \bar{A}_t, y\notin \bar{A}_t\cup B$ then $\bar{A}_t$ goes to
$(\bar{A}_t\backslash x)\cup y$ at rate $p(x,y)$ according to the
exponential clock of $\{x,y\}$. If $x\in \bar{A}_t, y\notin
\bar{A}_t\cup B^c$ and the exponential clock for $\{x,y\}$ goes
off then $\bar{A}_t$ goes to either $\bar{A}_t\backslash x$ if
$\eta_t(x)=1$ or the cemetery state $\Delta$ if $\eta_t(x)=0$.
Since the values of $\eta_t(x)$ and $\eta_t(y)$ switch when the
clock for $\{x,y\}$ goes off, we will assume that the evaluation
of $\eta_t(x)$ is taken before the switch.

For a fixed $T>0$, we define another process $\bar{A}_t^T$ to
follow the evolution described above except that it runs backwards
in time from $T$ to $0$ while $\eta_s$ runs forward in time; when
the exponential clock for $\{x,y\}$ goes off, the evaluation of
$\eta_s(x)$ takes place after the switching of $\eta_s(x)$ and
$\eta_s(y)$ at time $s=T-t$ takes place.  Setting
$\eta(\Delta)\equiv 0$, we then have following analog of
(\ref{q12}) for the quasi-symmetric process $\eta_t$:
\begin{equation}\label{q16}
\{\eta_T(x)=1\text{ for all }x\in
\bar{A}_0^T\}=\{\eta_0(x)=1\text{ for all }x\in \bar{A}_T^T\}.
\end{equation}

The processes $A_t$ and $\bar{A}_t$ are coupled so that they start
from the set $A\in\mathcal{S}_k$ (where $A\cap B=\emptyset$) and
move together as much as possible; likewise for the processes
$A_t^T$ and $\bar{A}_t^T$. Therefore denote
$$\mathcal{N}_A=\{\bar{A}_t,\text{ starting from }A,\text{ equals
}A_t\text{ for all }t\ge 0\}$$ and
$$\mathcal{N}_A^T=\{\bar{A}_t^T,\text{ starting from }A,\text{ equals
}A_t^T\text{ for all }t\in[0,T]\}.$$

In order to simplify some of the notation in the proof below
define the functions
$$f_{\bar{A}_T}(\bar{\omega})=\left\{
\begin{array}{ll}
1 &\text{ if } \eta_0(x)=1 \text{ for all }x\in \bar{A}_T(\bar{\omega})\\
0 &\text{ otherwise}
\end{array}
\right.$$ and
$$f_{\bar{A}_T^T}(\bar{\omega})=\left\{
\begin{array}{ll}
1 &\text{ if } \eta_T(x)=1 \text{ for all }x\in \bar{A}_T^T(\bar{\omega})\\
0 &\text{ otherwise}
\end{array}
\right.$$ The functions $f_{A_T}(\omega)$ and $f_{A_T^T}(\omega)$
are defined similarly.

\begin{proof}[Proof of Theorem \ref{qtheorem2}]

We prove part (b) first. Choose $\mu\in\mathcal{I}$. Using the
duality of symmetric exclusion given in (\ref{q12}) we have
\begin{eqnarray*}
\mu\{\eta(x)=1\text{ for all }x\in A\}-P(\mathcal{N}_A^c)&\le&
\int f_{A_T}  1_{\mathcal{N}_A}\, dP_\mu.
\end{eqnarray*}
Since $A_t=\bar{A}_t$ on $\mathcal{N}_{A}$ we get that
\begin{equation*}
\int f_{A_T}  1_{\mathcal{N}_A}\, dP_\mu\le \int
1_{\{\eta_0(x)=1\forall x\in
\bar{A}_T^T(\bar{\omega})\}}1_{\mathcal{N}_A^T} \, d\bar{P}_{\mu}.
\end{equation*}
Notice that the left-hand side is constant in $T$ since
$\mu\in\mathcal{I}$. Taking the limsup of both sides and using
(\ref{q16}) we get that
\begin{eqnarray*}
\int f_{A_T}  1_{\mathcal{N}_A}\, dP_\mu &\le&
\limsup_{T\rightarrow\infty}\int 1_{\{\eta_0(x)=1\forall x\in
\bar{A}_T^T(\bar{\omega})\}}1_{\mathcal{N}_A^T} \,
d\bar{P}_{\mu}\\
&=& \limsup_{T\rightarrow\infty}\int 1_{\{\eta_T(x)=1\forall x\in
A\}}1_{\mathcal{N}_A^T} \, d\bar{P}_\mu.
\end{eqnarray*}
By Theorem \ref{qtheorem} (proven in Section \ref{3}) we have that
$\lim_{t\rightarrow\infty} \mu \bar{S}(t)=\bar{\mu}$ so that the
right-hand side is less than or equal to
$\bar{\mu}\{\eta(x)=1\text{ for all }x\in A\}$ (in actuality we do
not require Theorem \ref{qtheorem} here since a Cesaro limit works
just as well, but it certainly simplifies things). Combining the
above arguments gives us
$$\mu\{\eta(x)=1\text{ for all }x\in A\}-P(\mathcal{N}_A^c)\le
\bar{\mu}\{\eta(x)=1\text{ for all }x\in A\}.$$

Similarly we have that
\begin{eqnarray*}
\bar{\mu}\{\eta(x)=1\text{ for all }x\in
A\}-P(\mathcal{N}_A^c)&\le& \limsup_{T\rightarrow\infty}\int
1_{\{\eta_T(x)=1\forall x\in
\bar{A}_0^T(\bar{\omega})\}}1_{\mathcal{N}_A^T} \, d\bar{P}_\mu\\
&=& \limsup_{T\rightarrow\infty}\int 1_{\{\eta_0(x)=1\forall x\in
\bar{A}_T^T(\bar{\omega})\}}1_{\mathcal{N}_A^T} \, d\bar{P}_\mu\\
&\le&\lim_{T\rightarrow\infty}E^A\int 1_{\{\eta(x)=1\forall x\in A_T\}} \, d\mu\\
&=&\mu\{\eta(x)=1\text{ for all }x\in A\}
\end{eqnarray*}
altogether giving us
\begin{eqnarray*}\label{q14}
|\bar{\mu}\{\eta(x)=1\text{ for all }x\in A\}-\mu\{\eta(x)=1\text{
for all }x\in A\}|\le P(\mathcal{N}_A^c).
\end{eqnarray*}
We complete the proof of part (b) by noting that Proposition
\ref{qlem5} tells us $\lim_{n\rightarrow\infty}
P(\mathcal{N}_{A^{n}}^c)=0$ for all $k$ and all sequences
$\{A^n\}$, $A^n\in\mathcal{S}_k$ such that each $x\in\mathcal{S}$
is in finitely many $A^n$.

The proof of part (a) is similar. Pick
$\bar{\mu}\in\bar{\mathcal{I}}$. Since $A_t=\bar{A}_t$ on
$\mathcal{N}_{A}$ we have for all $t\ge 0$ that
\begin{equation*}\label{q13}
\int f_{A_t} \, dP_{\bar{\mu}}-P(\mathcal{N}_{A}^c) \le \int
f_{\bar{A}_t}1_{\mathcal{N}_{A}} \, d\bar{P}_{\bar{\mu}}.
\end{equation*}
By (\ref{q12}) we have
\begin{equation*}
\int f_{A_t} \, dP_{\bar{\mu}}=E^A\int 1_{\{\eta(x)=1\forall x\in
A_t\}} \, d\bar{\mu}= \int 1_{\{\eta(x)=1\forall x\in A\}} \,
d\bar{\mu} S(t)
\end{equation*}
where $S(t)$ is the semigroup of symmetric exclusion. Since
$\bar{\mu}\in\bar{\mathcal{I}}$ we also have that
\begin{eqnarray*}
\int f_{\bar{A}_T}1_{\mathcal{N}_{A}} \, d\bar{P}_{\bar{\mu}}
%&\le&\int f_{\bar{A}_T^T}1_{\mathcal{N}_{A}^T}\,
%d\bar{P}_{\bar{\mu}}\\
&\le& \int  1_{\{\eta_0(x)=1\forall x\in
\bar{A}_T^T(\bar{\omega})\}}\, d\bar{P}_{\bar{\mu}}\\
&=& \int  1_{\{\eta_T(x)=1\forall x\in
\bar{A}_0^T(\bar{\omega})\}}\, d\bar{P}_{\bar{\mu}}\\
&=&\bar{\mu}\{\eta(x)=1\text{ for all }x\in A\}
\end{eqnarray*}
where the first equality follows from (\ref{q16}). Altogether we
have that
\begin{equation}\label{q20}
\int 1_{\{\eta(x)=1\forall x\in A\}} \, d\bar{\mu}
S(t)-P(\mathcal{N}_{A}^c) \le \bar{\mu}\{\eta(x)=1\text{ for all
}x\in A\}.
\end{equation}

Choose a sequence of times $\{T_n\}$ going to infinity so that
$$\lim_{n\rightarrow\infty} \frac{1}{T_n}\int_0^{T_n}\bar{\mu} S(t)
dt$$ converges to some $\mu\in\mathcal{I}$
%. By using the sequence
%$\{T_n\}$ to take the Cesaro limit of the left-hand side of
so that (\ref{q20}) gives us
\begin{equation*}\label{q9}
\mu\{\eta(x)=1\text{ for all }x\in A\}-P(\mathcal{N}_A^c)\le
\bar{\mu}\{\eta(x)=1\text{ for all }x\in A\}.
\end{equation*}

Similarly
\begin{eqnarray*}
\bar{\mu}\{\eta(x)=1\text{ for all }x\in
A\}-P(\mathcal{N}_A^c)&\le&
\lim_{n\rightarrow\infty}\frac{1}{T_n}\int_0^{T_n}\int
1_{\{\eta_T(x)=1\forall x\in\bar{A}_0^T(\bar{\omega})\}}
1_{\mathcal{N}_A^T} \, d\bar{P}_{\bar{\mu}}\, dT\\
&=& \lim_{n\rightarrow\infty}\frac{1}{T_n}\int_0^{T_n}\int
1_{\{\eta_0(x)=1\forall
x\in\bar{A}_T^T(\bar{\omega})\}}1_{\mathcal{N}_A^T} \,
d\bar{P}_{\bar{\mu}}\, dT
\\
&\le& \lim_{n\rightarrow\infty}\frac{1}{T_n}\int_0^{T_n}\int
1_{\{\eta_0(x)=1\forall x\in
{A}_T^T({\omega})\}} \, d{P}_{\bar{\mu}}\, dT\\
&=& \lim_{n\rightarrow\infty}\frac{1}{T_n}\int_0^{T_n}\int
1_{\{\eta_T(x)=1\forall x\in
{A}_0^T({\omega})\}} \, d{P}_{\bar{\mu}}\, dT\\
&=&\mu\{\eta(x)=1\text{ for all }x\in A\}
\end{eqnarray*}
giving us
\begin{eqnarray*}
|\bar{\mu}\{\eta(x)=1\text{ for all }x\in A\}-\mu\{\eta(x)=1\text{
for all }x\in A\}|\le P(\mathcal{N}_A^c).
\end{eqnarray*}
As in the proof of part (a) of the theorem, we apply Proposition
\ref{qlem5} to get that (\ref{q11}) holds.
\end{proof}

\section{The Infinitesimal Coupling}\label{im}
The main tool used in the proof of the Theorems \ref{qtheorem} and
\ref{newtheorem} is the so called \textit{infinitesimal coupling}
of the process $\eta_t$. In this section we will describe the
infinitesimal coupling and present some results concerning this
coupling.

The infinitesimal coupling of the process $\eta_t$ follows the
motion of the \textit{basic coupling} (defined below) for the two
processes $\eta_t$ and $\xi_t^s$ having joint initial measure
$\tilde{\nu}$ (also defined below). The marginal process $\xi_t^s$
can be thought of as an approximation of $\eta_{t+s}$ for small
values of $s$.

Let us now define the basic coupling of two exclusion processes
$\eta_t$ and $\xi_t$ having the same generator. Simply put, the
basic coupling is the coupling which allows $\eta_t$ and $\xi_t$
to move together as much as possible. The generator for the basic
coupling is the closure of the operator $\tilde{\Omega}$ defined
on $\mathcal{D}(X\times X)$:
\begin{eqnarray*}
&&\tilde{\Omega}f(\eta,\xi)=\sum_{\eta(x)=\xi(x)=1,\eta(y)=\xi(y)=0}{\bar{p}(x,y)[f(\eta_{xy},\xi_{xy})-f(\eta,\xi)]}\\
&+&\sum_{\eta(x)=1,\eta(y)=0 \text{ and } (\xi(y)=1 \text{ or } \xi(x)=0)}{\bar{p}(x,y)[f(\eta_{xy},\xi)-f(\eta,\xi)]}\\
&+&\sum_{\xi(x)=1,\xi(y)=0 \text{ and } (\eta(y)=1 \text{ or }
\eta(x)=0)}{\bar{p}(x,y)[f(\eta,\xi_{xy})-f(\eta,\xi)]}.
\end{eqnarray*}

The initial measure $\tilde{\nu}$ depends on the transition kernel
of the process.  To describe $\tilde{\nu}$, we will consider the
following simple kernel: Start with an irreducible transition
kernel $p(x,y)$ on $\mathcal{S}$ (not necessarily symmetric). Pick
an ordered pair $(u,v)$. Choosing $\epsilon>0$, we can define
$\bar{p}(x,y)$ by
\begin{equation}\label{q1}
\bar{p}(u,v)=p(u,v)+\epsilon,\ \  \bar{p}(x,y)=p(x,y) \text{
elsewise}.
\end{equation}

In order to simplify the description of $\tilde{\nu}$, we will
assume throughout most of this section that our transition kernel
is given by (\ref{q1}). It is under this assumption that we will
explicitly describe $\tilde{\nu}$ and prove the lemmas.  At the
end of the section we will give an argument that extends the
results to general perturbed kernels.

We are ready to describe $\tilde{\nu}$ under the assumption of
(\ref{q1}). Following Andjel, Bramson, and Liggett(1988), the
basic idea is to couple a given measure $\mu\in\mathcal{I}$
together with $\mu \Bar{S}(s)$ for small values of $s$ (in
particular, we impose the restriction $s<\frac{1}{\epsilon}$). The
problem is that one cannot explicitly write out the distribution
of $\mu \Bar{S}(s)$; however, it turns out that a first order
approximation to $\mu \Bar{S}(s)$ is good enough. Therefore, we
think of $\mu^s$ as some measure $\mu \Bar{S}(s)+o(s)$ as
$s\rightarrow 0$. Throughout the rest of the section $\mu$ will be
the marginal distribution of $\tilde{\nu}$ corresponding to
$\eta_0$ and $\mu^s$ will be the marginal distribution of
$\tilde{\nu}$ corresponding to $\xi^s_0$.

The measures $\mu^s$ and $\tilde{\nu}$ will be defined in such a
way that $\tilde{\nu}$ has a small number of discrepancies (a
\textit{discrepancy} occurs when $\eta(x)\neq \xi^s(x)$). This is
because the idea is to let the coupled process run according to
the basic coupling and analyze the behavior of the discrepancies.
In fact, it is by analyzing the behavior of the discrepancies that
we will be able to prove that the measure
$\lim_{t\rightarrow\infty} \mu \Bar{S}(t)$ exists for all
$\mu\in\mathcal{I}$.

Let us now explicitly describe $\mu^s$.  If $D$ is the set
$\{\eta_{0}(u)=1,\eta_{0}(v)=0\}$ then define $\mu_D$ and
$\mu_{D^c}$ by conditioning $\mu$ on the events $D$ and $D^c$.
Also, define $\hat{\mu}_D$ to be the measure that is exactly
$\mu_D$ except that $\xi^s_{0}(u)=0$ and $\xi^s_{0}(v)=1$. We then
have
$$\mu^s=[\mu\{D\}(1-s\epsilon)]\mu_D+[\mu\{D\}s\epsilon]\hat{\mu}_D+[\mu\{D^c\}]\mu_{D^c}.$$
Note that this measure is well-defined for $s<\frac{1}{\epsilon}$.

Let ${\mu}_D$ and $\hat{\mu}_D$ be coupled in such a way that they
agree everywhere except at $u$ and $v$.  The coupling measure
$\tilde{\nu}$ is just the coupling of $\eta_0$ and $\xi_0^s$ such
that the two marginals agree everywhere except on a set of measure
$\mu\{D\}s\epsilon$ where we use the coupling of ${\mu}_D$ and
$\hat{\mu}_D$ described in the previous sentence.  In particular,
the distribution for
\begin{equation}\label{margin}
\left(
\begin{array}{ccc}
\xi^s_{0}(u)&\xi^s_{0}(v)\\
\eta_0(u)&\eta_0(v)
\end{array}
\right)
\end{equation}
is given by

$
\begin{array}{ll}
\;\text{ Value}&\text{Probability}\\
\left(
\begin{array}{ccc}
1&1\\
1&1
\end{array}
\right) &\mu\{\eta_{0}(u)=1,\eta_{0}(v)=1\}\\
\\
\left(
\begin{array}{ccc}
1&0\\
1&0
\end{array}
\right) &\mu\{D\}(1-s\epsilon)\\
\\
\left(
\begin{array}{ccc}
0&1\\
0&1
\end{array}
\right) &\mu\{\eta_{0}(u)=0,\eta_{0}(v)=1\}\\
\\
\left(
\begin{array}{ccc}
0&0\\
0&0
\end{array}
\right) &\mu\{\eta_{0}(u)=0,\eta_{0}(v)=0\}\\
\\
\left(
\begin{array}{ccc}
0&1\\
1&0
\end{array}
\right) &\mu\{D\}s\epsilon.
\end{array}
$

As desired, up to first order in $s$,
$(\xi^s_{0}(u),\xi^s_{0}(v))$ has the same distribution as
$(\eta_s(u),\eta_s(v))$ under $\mu$. This is what lies behind the
next lemma.
\begin{lemma}\label{qlem1}
Suppose $\mu\in\mathcal{I}$.  Then for any
$f\in\mathcal{D}(\{0,1\}^\mathcal{S})$,
\begin{equation*}
\lim_{s\rightarrow 0} \frac{Ef(\xi^s_{0})-\int f\, d\mu
\Bar{S}(s)}{s}=0.
\end{equation*}
\end{lemma}

\begin{proof}
Let $\Omega$ be the generator with respect to $p(x,y)$ and
$\bar{\Omega}$ be the generator with respect to $\bar{p}(x,y)$.
Using (\ref{qgenerator}), (\ref{q1}) and the fact that
$\mu\in\mathcal{I}$ we have
\begin{eqnarray*}
\int \bar{\Omega} f \,d\mu&=& \int
\sum_{x,y}{\bar{p}(x,y)\xi(x)(1-\xi(y))[f(\xi_{xy})-f(\xi)]}\,
d\mu\\
&=&\int\sum_{x,y}{{p}(x,y)\xi(x)(1-\xi(y))[f(\xi_{xy})-f(\xi)]}\,
d\mu +\int \epsilon\,\xi(u)(1-\xi(v))[f(\xi_{uv})-f(\xi)]\,
d\mu\\
&=&\int\Omega f\, d\mu +\int
\epsilon\,\xi(u)(1-\xi(v))[f(\xi_{uv})-f(\xi)]\,
d\mu\\
&=&\int\epsilon\,\xi(u)(1-\xi(v))[f(\xi_{uv})-f(\xi)] \, d\mu.
\end{eqnarray*}

But now, using the explicit expression for the distribution of
$\xi^s_0$, we also get for $s>0$ that
\begin{equation*}
\frac{Ef(\xi^s_{0})-\int f\, d\mu }{s}=\int
\epsilon\,\xi(u)(1-\xi(v))[f(\xi_{uv})-f(\xi)] \, d\mu=\int
\bar{\Omega} f \,d\mu.
\end{equation*}

By the definition of the generator
\begin{equation*}
\int \bar{\Omega} f \,d\mu=\lim_{s\rightarrow 0} \frac{\int f\,
d\mu \Bar{S}(s)-\int f\, d\mu}{s}.
\end{equation*}

Combining the last two equations gives us
\begin{equation*}
\lim_{s\rightarrow 0} \frac{Ef(\xi^s_{0})-\int f\, d\mu
\Bar{S}(s)}{s}=0.
\end{equation*}
\end{proof}
%It is important to note that the above lemma does not need that
%the original unperturbed transition kernel $p(x,y)$ is symmetric.
%In the generalization of this lemma, we only require that
%$\bar{p}(x,y)$ differs from $p(x,y)$ for finitely many pairs
%$\{(x_i,y_i)\}$ and that $\mu$ is invariant for the exclusion
%process with respect to the kernel $p(x,y)$.

Define $(\hat{\eta}_{t},\hat{\xi}^s_{t})$ by conditioning
$(\eta_t,\xi^s_t)$ on the event that
\begin{equation*}
\left(
\begin{array}{cc}
\xi^s_0(u)&\xi_0^s(v)\\
\eta_0(u)&\eta_0(v)
\end{array}
\right)=\left(
\begin{array}{ccc}
0&1\\
1&0
\end{array}
\right).
\end{equation*}
This is the only event for which $\eta_0$ and $\xi^s_0$ differ.
Note that after conditioning, the distribution of the coupling no
longer depends on $s$.

The proof of the next lemma follows that of Lemma 3.4 in Andjel,
Bramson, and Liggett(1988).

\begin{lemma}\label{coupling2}
If $A$ is any finite subset of $\mathcal{S}$ then
\begin{equation*}
 \frac{d}{dt}\mu
\Bar{S}(t)\{\eta:\eta(x)=1\text{ for all }x\in
A\}={\epsilon\mu\{D\}} E [\prod_{x\in A}
\hat{\xi}_t^s(x)-\prod_{x\in A}\hat{\eta}_t(x)]
\end{equation*}
\end{lemma}

\begin{proof}
Let
\begin{equation*}
f_A(\eta)=\prod_{x\in A} \eta(x)=\left\{
\begin{array}{ll}
1 &\text{ if } \eta(x)=1 \text{ for all }x\in {A}\\
0 &\text{ otherwise.}
\end{array}
\right.
\end{equation*}
Then $f_A\in \mathcal{D}({X})$, so $f_A^t=\Bar{S}(t)f_A$ is also
in $\mathcal{D}({X})$ by Theorem I.3.9 of IPS. Letting $\mu^t=\mu
\bar{S}(t)$, we compute
\begin{eqnarray*}
&&\frac{d}{dt}\mu \Bar{S}(t)\{\eta(x)=1\text{ for all }x\in
A\}\\
&=&\lim_{s\rightarrow 0}\frac{1}{s}[\mu^{t+s}
\{\eta(x)=1\text{ for all }x\in A\}-\mu^t\{\eta(x)=1\text{ for all }x\in A\}]\\
&=&\lim_{s\rightarrow 0}\frac{1}{s}[\int f_A\,
d\mu^{t+s}-\int f_A\,d\mu^t]\\
&=&\lim_{s\rightarrow 0}\frac{1}{s}[\int f_A^t\,
d\mu^{s}-\int f_A^t\,d\mu]\\
&=&\lim_{s\rightarrow 0}\frac{E f_A^t(\xi_0^s)-\int
f_A^t\,d\mu}{s}
\end{eqnarray*}
where the last equality follows from Lemma \ref{qlem1}.  This in
turn equals
\begin{eqnarray*}
\lim_{s\rightarrow 0}\frac{E f_A^t(\xi_0^s)-E
f_A^t(\eta_0)}{s}&=&\lim_{s\rightarrow 0}\frac{E
f_A(\xi_t^s)-E f_A(\eta_t)}{s}\\
&=&\lim_{s\rightarrow 0}\frac{1}{s} E [\prod_{x\in A}
\xi_t^s(x)-\prod_{x\in A}\eta_t(x)]\\
&=&{\epsilon\mu\{D\}} E
[\prod_{x\in A} \hat{\xi}_t^s(x)-\prod_{x\in A}\hat{\eta}_t(x)]
\end{eqnarray*}
\end{proof}

Let $({\eta}_{t}^{(z)},{\xi}_t^{(z)})$ be a process that runs
according to the basic coupling for $z=u,v$.  Its initial
distribution is such that both the marginal distributions
(corresponding to ${\eta}_{0}^{(z)}$ and ${\xi}_0^{(z)}$) are
equal to the measure $\mu_D$ except that we force
$\xi_0^{(z)}(z)=1, \eta_0^{(z)}(z)=0$. As usual, the initial
distribution is coupled such that $\xi_0^{(z)}(x)=\eta_0^{(z)}(x)$
for all $x\neq z$.
\begin{corollary}\label{coupling4}
If $A$ is any finite subset of $\mathcal{S}$ then
\begin{equation*}
 |\frac{d}{dt}\mu
\Bar{S}(t)\{\eta:\eta(x)=1\text{ for all }x\in A\}|\le
\epsilon\mu\{D\}\sum_{z=0,1}\sum_{x\in
A}E[{\xi}^{(z)}_t(x)-{\eta}^{(z)}_t(x)].
\end{equation*}
\end{corollary}

\begin{proof}
\begin{eqnarray*}
{\epsilon\mu\{D\}} |E [\prod_{x\in A} \hat{\xi}_t^s(x)-\prod_{x\in
A}\hat{\eta}_t(x)]|&\le& {\epsilon\mu\{D\}} E |\prod_{x\in A}
\hat{\xi}_t^s(x)-\prod_{x\in A}\hat{\eta}_t(x)|\\
 &\le&\epsilon\mu\{D\}\sum_{x\in A}
P(\hat{\xi}_t(x)\neq\hat{\eta}_t(x))
\\
&\le&\epsilon\mu\{D\}\sum_{z=0,1}\sum_{x\in A}
E({\xi}^{(z)}_t(x)-{\eta}^{(z)}_t(x)).
\end{eqnarray*}
The last inequality is due to a property given by the basic
coupling: when the two discrepancies
\begin{equation*}
\left(\begin{array}{c}
\xi_T^s(x)=1\\
\eta_T(x)=0
\end{array}\right)
\text{ and } \left(\begin{array}{c}
\xi_T^s(x)=0\\
\eta_T(x)=1
\end{array}\right)
\end{equation*}
meet, they cancel each other out to result in no discrepancies for
all $t\ge T$.
\end{proof}
We now give an argument that extends the infinitesimal coupling
and the results to a general perturbed kernel.  The first thing is
to realize that if $\epsilon$ is negative, we can obtain analogs
of the two lemmas if we make the following changes to the
distribution of (\ref{margin}):

$
\begin{array}{ll}
\;\text{ Value}&\text{Probability}\\
\left(
\begin{array}{ccc}
1&0\\
1&0
\end{array}
\right) &\mu\{D\}\\
\\
\left(
\begin{array}{ccc}
0&1\\
0&1
\end{array}
\right) &\mu\{\eta_{0}(u)=0,\eta_{0}(v)=1\}-\mu\{D\}s|\epsilon|\\
\\
\left(
\begin{array}{ccc}
1&0\\
0&1
\end{array}
\right) &\mu\{D\}s|\epsilon|\\
\\
\left(
\begin{array}{ccc}
0&1\\
1&0
\end{array}
\right) &0.
\end{array}
$
\newline
Here we impose the restriction
$s<\frac{\mu\{\eta_{0}(u)=0,\eta_{0}(v)=1\}}{\mu\{D\}|\epsilon|}$.

Next we see that if there are multiple differences between
$p(x,y)$ and $\bar{p}(x,y)$, we can superimpose the changes to the
distribution of $\tilde{\nu}$ to get analogs of the two lemmas.
For instance if
\begin{equation*}
\bar{p}(w,y)=p(w,y)+\epsilon_1\text{ and
}\bar{p}(w,z)=p(w,z)+\epsilon_2\text{ where }\epsilon_i>0,
\end{equation*}
then when $s<\frac{1}{\epsilon_1+\epsilon_2}$, the distribution of
the coupling at $(w,y,z)$ at time $0$ is identical to the marginal
measures for $(\eta_0(w),\eta_0(y),\eta_0(z))$ and for
$(\xi_0^s(w),\xi_0^s(y),\xi_0^s(z))$, except at the values in the
table below:

$
\begin{array}{ll}
\;\text{ Value}&\text{Probability}\\
\left(
\begin{array}{ccc}
1&0&0\\
1&0&0
\end{array}
\right) &\mu\{\eta_0(w)=1,\eta_0(y)=0,\eta_0(z)=0\}[1-s(\epsilon_1+\epsilon_2)]\\
\\
\left(
\begin{array}{ccc}
1&0&1\\
1&0&1
\end{array}
\right) &\mu\{\eta_0(w)=1,\eta_0(y)=0,\eta_0(z)=1\}(1-s\epsilon_1)\\
\\
\left(
\begin{array}{ccc}
1&1&0\\
1&1&0
\end{array}
\right) &\mu\{\eta_0(w)=1,\eta_0(y)=1,\eta_0(z)=0\}(1-s\epsilon_2)\\
\\
\left(
\begin{array}{ccc}
0&1&0\\
1&0&0
\end{array}
\right) &\mu\{\eta_0(w)=1,\eta_0(y)=0,\eta_0(z)=0\}s\epsilon_1\\
\\
\left(
\begin{array}{ccc}
0&1&1\\
1&0&1
\end{array}
\right) &\mu\{\eta_0(w)=1,\eta_0(y)=0,\eta_0(z)=1\}s\epsilon_1\\
\\
\left(
\begin{array}{ccc}
0&0&1\\
1&0&0
\end{array}
\right) &\mu\{\eta_0(w)=1,\eta_0(y)=0,\eta_0(z)=0\}s\epsilon_2\\
\\
\left(
\begin{array}{ccc}
0&1&1\\
1&1&0
\end{array}
\right) &\mu\{\eta_0(w)=1,\eta_0(y)=1,\eta_0(z)=0\}s\epsilon_2.
\end{array}
$

%If there are $n$ perturbations $\{\epsilon_1, \ldots,
%\epsilon_n\}$ then we have $n$ events $\{E_1,\ldots,E_n\}$ for
%which $\eta_0$ and $\xi^s_0$ differ. Defining
%$({\eta}_{t}^i,{\xi}^i_{t})$ by conditioning $(\eta_t,\xi^s_t)$ on
%the event $E_i$ we can extend Lemma \ref{coupling2}:

%\begin{lemma}
%If $A$ is any finite subset of $\mathcal{S}$ then there exist
%constants $C_i$ such that
%\begin{equation*}
% \frac{d}{dt}\mu
%\Bar{S}(t)\{\eta:\eta(x)=1\text{ for all }x\in A\}=\sum_{i=1}^n
%C_i E [\prod_{x\in A} {\xi}_t^i(x)-\prod_{x\in A}{\eta}^i_t(x)]
%\end{equation*}
%\end{lemma}

Recall that
\begin{equation*}
B=\{x\in\mathcal{S}: \bar{p}(x,y)\neq {p}(x,y) \text{ or
}\bar{p}(y,x)\neq {p}(y,x)\text{ for some }y\in\mathcal{S}\}.
\end{equation*}
If we define $({\eta}_{t}^{(z)},{\xi}_t^{(z)})$ for all $z\in B$
similarly to our previous definition, then we get the following
extension of Corollary \ref{coupling4}:
\begin{corollary}\label{coupling3}
If $A$ is any finite subset of $\mathcal{S}$ then there exists a
constant $C$ such that
\begin{eqnarray*}
|\frac{d}{dt}\mu \Bar{S}(t)\{\eta:\eta(x)=1\text{ for all }x\in
A\}| \le C\sum_{z\in B}\sum_{x\in
A}E[{\xi}^{(z)}_t(x)-{\eta}^{(z)}_t(x)].
\end{eqnarray*}
\end{corollary}
The proof of the corollary is essentially the same as that of
Corollary \ref{coupling4} so we only make the following remark. It
is important to note that a pair of discrepancies of opposite type
$\left(
\begin{array}{c}
1\\
0
\end{array}
\right)$ and $\left(
\begin{array}{c}
0\\
1
\end{array}
\right)$ occur together, but any two pairs do not occur at the
same time.  Therefore, we still have that the only interaction
between discrepancies is when two discrepancies of opposite type
cancel each other out.

\section{Proofs of Theorems \ref{qtheorem}-\ref{newtheorem2}}\label{3}
Assume throughout this section that $p(x,y)$ is a symmetric kernel
and $\bar{p}(x,y)$ is a corresponding quasi-symmetric kernel.
Also, let $Y_t$ be the continuous-time Markov chain with kernel
$p(x,y)$.

Given the process $({\eta}_{t}^{(z)},{\xi}_t^{(z)})$ described in
the previous section, let $Y_t^*$ mark the position at time $t$ of
the discrepancy that starts at $z$.  Notice that while the process
$Y_t^*$ is not a Markov process, the joint process
$(Y_t^*,\eta_t)$ is a Markov process. Let
\begin{equation*}
G^*(z,x)= E^{z}\int_0^\infty P(Y_t^*=x)\, dt
\end{equation*}
be the expected time that the discrepancy starting at $z$ spends
at $x$. We note here that the above expectation is taken over the
path space of the joint process $(Y_t^*,\eta_t)$ where the initial
distribution is taken to be the initial measure described
immediately following the proof of Lemma \ref{coupling2} (we have
used the notation $E^z$ to indicate that the discrepancy starts
from $z$). If $Y_n^*$ is the embedded discrete-time process for
$Y_t^*$, define
$$H^*(z,x)=\sup_\eta P^{(z,\eta)}(Y_n^*=x\text{ for some }n\ge
1).$$

\begin{lemma}\label{qlem4}
If $Y_t$ is transient then $G^*(z,x)<\infty$ for all
$z,x\in\mathcal{S}$.
\end{lemma}
\begin{proof}
If the discrepancy is at site $x$, it goes to $y$ at rate
$\bar{p}(x,y)$ when $\xi^{(z)}(y)=\eta^{(z)}(y)=0$ and at rate
$\bar{p}(y,x)$ when $\xi^{(z)}(y)=\eta^{(z)}(y)=1$.  But when
$x\notin B$, $\bar{p}(x,y)=\bar{p}(y,x)$ since $p(x,y)$ is
symmetric. Therefore when $Y_t^*\notin B$, $Y_t^*$ moves according
to the same transition rates as $Y_t$.

Couple $Y_t^*$ with the process $Y_t$ starting from $z$ so that
they move together as much as possible and let
$$E=\{\omega:Y_t^*(\omega)=Y_t(\omega) \text{ for all }t\ge 0, Y_n^*\neq z\text{ for all }n\ge 1\}.$$ Since $B$ is finite and $Y_t$ is transient, and since
$\bar{p}(x,y)>0$ whenever $p(x,y)>0$ (by the definition of
$\bar{p}(x,y)$), we have that $\inf_\eta P^{(z,\eta)}(E)>0$.

For each $x$ we have
\begin{eqnarray*}
H^*(x,x)&=&\sup_\eta P^{(x,\eta)}\left[\{Y_n^*=x\text{ for some
}n\ge 1\} \cap (E\cup
E^c)\right]\\
&=&\sup_\eta P^{(x,\eta)}(\{Y_n^*=x\text{ for some }n\ge 1\} \cap
E^c)\le 1-\inf_\eta P^{(x,\eta)}(E).
\end{eqnarray*}
Using the proof of Proposition 4-20 in Kemeny, Snell, and
Knapp(1976) we get that for some constant $C$,
$$G^*(z,x)\le C\sum_{k\ge
0}(H^*(x,x))^k<\infty.$$
\end{proof}

\begin{proof}[Proof of Theorem \ref{qtheorem}]
By the Inclusion-Exclusion Principle we need only show that for
each finite set $A\subset \mathcal{S}$,
\begin{equation}\label{q3}
\lim_{t\rightarrow\infty} \mu \Bar{S}(t) \{\eta:\eta(x)=1\text{
for all }x\in A\}
\end{equation}
exists.

Suppose to the contrary that there exists some $A$ for which
(\ref{q3}) does not exist. Then there exists a sequence $\{t_n\}$
going to infinity such that the set
\begin{equation*}
\{\mu \Bar{S}(t_n) \{\eta(x)=1\text{ for all }x\in A\}\}
\end{equation*}
has at least two different limit points. Therefore it must be that
\begin{equation*}
\int_0^\infty |\frac{d}{dt}\mu \Bar{S}(t)\{\eta(x)=1\text{ for all
}x\in A\}| dt=\infty.
\end{equation*}
On the other hand, by Corollary \ref{coupling3} and Lemma
\ref{qlem4},
\begin{eqnarray*}
\int_0^\infty |\frac{d}{dt}\mu \Bar{S}(t)\{\eta:\eta(x)=1\text{
for all }x\in A\}| dt&\le& C \int_0^\infty \sum_{z\in B}\sum_{x\in
A}E[{\xi}^{(z)}_t(x)-{\eta}^{(z)}_t(x)] dt\\
&\le& C\sum_{z\in B}\sum_{x\in A} G^*(z,x)<\infty,\nonumber
\end{eqnarray*}
a contradiction.  Therefore (\ref{q3}) exists for all finite $A$.
\end{proof}

\begin{proof}[Proof of Theorem \ref{newtheorem}]
We will prove only the case where $\bar{p}(u,v)=p(u,v)+\epsilon$
where $-p(u,v)<\epsilon<0$.  The proof of the case $\epsilon>0$ is
similar and left to the reader.

By the remarks following Theorem \ref{qtheorem2} we need only show
that for each $\mu\in\mathcal{I}$, $\lim_{t\rightarrow\infty}\mu
S(t)=\bar{\mu}$ is the only measure in $\bar{\mathcal{I}}$
satisfying the asymptotics given in (\ref{q11}). Suppose to the
contrary that there exists another measure
$\bar{\nu}\neq\bar{\mu}, \bar{\nu}\in\bar{\mathcal{I}}$ so that
both $\bar{\nu}$ and $\bar{\mu}$ satisfy (\ref{q11}) for the same
$\mu\in\mathcal{I}$.

As in Section \ref{im}, let
$$D=\{\eta:\eta(u)=1,
\eta(v)=0\}.$$ If we reverse the roles of $p(x,y)$ and
$\bar{p}(x,y)$ by thinking of $p(x,y)$ as a perturbation of
$\bar{p}(x,y)$ where $p(u,v)=\bar{p}(u,v)+|\epsilon|$ then all of
the results in Section \ref{im} are still true. Let
$(\bar{\xi}_t^s, \bar{\eta}_t)$ and
$(\tilde{\xi}_t^s,\tilde{\eta}_t)$ be processes similar to
$(\hat{\xi}_t^s, \hat{\eta}_t)$ except that they are defined with
respect to the measures $\bar{\mu}$ and $\bar{\nu}$ respectively.
By Lemma \ref{coupling2} and the Inclusion-Exclusion Principle, we
then have
\begin{eqnarray}\label{hello} \mu\{D\}-\bar{\mu}\{D\}&=&\int_0^\infty
\frac{d}{dt}\bar{\mu}
{S}(t)\{D\} dt\\
&=& {\epsilon\bar{\mu}\{D\}} \int_0^\infty[P(\bar{\xi}_t^s(u)=1,
\bar{\xi}_t^s(v)=0)- P(\bar{\eta}_t(u)=1, \bar{\eta}_t(v)=0)]
dt\nonumber
\end{eqnarray}
and
\begin{eqnarray*}
\mu\{D\}-\bar{\nu}\{D\}&=&\int_0^\infty \frac{d}{dt}\bar{\nu}
{S}(t)\{D\} dt\\
&=& {\epsilon\bar{\nu}\{D\}}\int_0^\infty [P(\tilde{\xi}_t^s(u)=1,
\tilde{\xi}_t^s(v)=0)- P(\tilde{\eta}_t(u)=1,
\tilde{\eta}_t(v)=0)] dt
\end{eqnarray*}
Note that convergence of the integrals above follows from the
transience of the discrepancies $\left(
\begin{array}{c}
1\\
0
\end{array}
\right)$ and $\left(
\begin{array}{c}
0\\
1
\end{array}
\right)$ in the processes $(\bar{\xi}_t^s, \bar{\eta}_t)$ and
$(\tilde{\xi}_t^s,\tilde{\eta}_t)$.

Since the semigroup $S(t)$ acts with respect to the symmetric
kernel $p(x,y)$, it follows that the discrepancies in the process
$(\bar{\xi}_t^s, \bar{\eta}_t)$ (respectively
$(\tilde{\xi}_t^s,\tilde{\eta}_t)$) behave exactly like Markov
chains with kernels $p(x,y)$. In particular this means that we can
couple the discrepancy $\left(
\begin{array}{c}
1\\
0
\end{array}
\right)$ in the process $(\bar{\xi}_t^s, \bar{\eta}_t)$ with the
discrepancy $\left(
\begin{array}{c}
1\\
0
\end{array}
\right)$ in the process $(\tilde{\xi}_t^s,\tilde{\eta}_t)$ so that
they always move together (similarly for the $\left(
\begin{array}{c}
0\\
1
\end{array}
\right)$'s).

Applying this coupling to (\ref{hello}), we get that there exists
a constant $K$ such that
\begin{eqnarray}\label{goodbye} \mu\{D\}-\bar{\mu}\{D\}= K{\bar{\mu}\{D\}}
\,\text{ and }\, \mu\{D\}-\bar{\nu}\{D\}=K\bar{\nu}\{D\}.
\end{eqnarray}
However, these equations can only hold when either
$\bar{\mu}\{D\}=\bar{\nu}\{D\}$ or when $\mu\{D\}=0$.  In the
latter case, irreducibility implies that $\mu$ must be either
$\delta_0$ or $\delta_1$ (the measures that concentrate on
$\{\eta:\eta(x)\equiv 0\}$ and $\{\eta:\eta(x)\equiv 1\}$), but in
these cases it is clear that $\bar{\mu}$ and $\bar{\nu}$ must also
be equal to either $\delta_0$ or $\delta_1$ completing the proof.
\end{proof}
We comment that the proof of an extension of Theorem
\ref{newtheorem} to the general case breaks down when going from
(\ref{hello}) to (\ref{goodbye}) since the extra variables
introduced prevent us from getting a ``unique solution''.

\begin{proof}[Proof of Theorem \ref{newtheorem2}]
Pick $\bar{\mu}\in\bar{\mathcal{I}}$ and choose a sequence
$\{T_n\}$ such that
\begin{equation*}
\lim_{n\rightarrow\infty}\frac{1}{T_n}\int_0^{T_n} \bar{\mu}
{S}(t) dt=\mu\in\mathcal{I}.
\end{equation*}

If there are $n$ perturbations
$\{\bar{p}(x_1,y_1)=p(x_1,y_1)+\epsilon_1, \ldots,
\bar{p}(x_n,y_n)=p(x_n,y_n)+\epsilon_n\}$ then there are $n$
events $\{E_1,\ldots,E_n\}$ for which $\eta_0$ and $\xi^s_0$
differ. Define $({\eta}_{t}^i,{\xi}^i_{t})$ by conditioning
$(\eta_t,\xi^s_t)$ on the event $E_i$. The proof of Lemma
\ref{coupling2} can be generalized to show that for any finite
$A\subset\mathcal{S}$ there exist constants $C_i$ such that
\begin{equation*}
 \frac{d}{dt}\bar{\mu}
{S}(t)\{\eta:\eta(x)=1\text{ for all }x\in A\}=\sum_{i=1}^n C_i E
[\prod_{x\in A} {\xi}_t^i(x)-\prod_{x\in A}{\eta}^i_t(x)]
\end{equation*}
giving us
\begin{eqnarray*}\label{abc}
\mu\{\eta=1\text{ on }A\}-\bar{\mu}\{\eta=1\text{ on
}A\}&=&\lim_{n\rightarrow\infty}\frac{1}{T_n}\int_0^{T_n}\int_0^s
\sum_{i=1}^n C_i E [\prod_{x\in A} {\xi}_t^i(x)-\prod_{x\in
A}{\eta}^i_t(x)] dt \, ds
\end{eqnarray*}

As in the proof of Theorem \ref{newtheorem}, the semigroup $S(t)$
acts with respect to the symmetric kernel $p(x,y)$ so the
discrepancies behave exactly like Markov chains with kernels
$p(x,y)$. The expected amount of time that the discrepancies
between ${\xi}_t^i$ and ${\eta}^i_t$ spend at a given
$x\in\mathcal{S}$ up to time $T$ is then given by the Green's
function $G_T(z,x)$ where $z$ is the position of the discrepancy
at $t=0$. Since discrepancies come in pairs of opposite type and
since $$a(x-z)-a(x-y)=\lim_{T\rightarrow\infty}
[G_T(y,x)-G_T(z,x)],$$ we have that
$$|\mu\{\eta=1\text{ on }A\}-\bar{\mu}\{\eta=1\text{ on
}A\}|\le\sum_{i=1}^n\sum_{x\in A} |C_i(a(x-x_i)-a(x-y_i))|.$$
Fixing $k$ and $\epsilon>0$ we can choose $N$ large enough so that
for all $A\in\mathcal{S}_k$ with $d(A,o)>N$, the right-hand side
in the above inequality is less than $\epsilon$. We therefore have
that (\ref{q11.1}) is satisfied.
\end{proof}

\textbf{Acknowledgement}. The author thanks his advisor, Thomas M.
Liggett, for his continual support and for the many discussions
that led to the writing of this paper.

\end{document}